\documentclass[11pt,reqno]{amsart}
\usepackage{amssymb}
\usepackage{eurosym}
\usepackage{amsfonts}
\usepackage{amsmath}
\usepackage{graphicx}
\usepackage{amscd}
\usepackage{fullpage}
\usepackage{amsfonts}
\usepackage{mathrsfs,float}
\restylefloat{figure}

\theoremstyle{plain}

\numberwithin{equation}{section}

\begin{document}
\title[Oscillation of Equations with Non-monotone Arguments]{ITERATIVE
OSCILLATION TESTS FOR\ DIFFERENCE\\
EQUATIONS WITH SEVERAL NON-MONOTONE ARGUMENTS
}
\author{E. BRAVERMAN$^1$}
\address{Department of Mathematics and Statistics\\
University of Calgary\\
2500 University Drive N. W., Calgary, Canada AB T2N 1N4}
\email{maelena@math.ucalgary.ca}
\author{G. E. CHATZARAKIS}
\address{Department of Electrical and Electronic Engineering Educators\\
School of Pedagogical and Technological Education (ASPETE)\\
14121, N. Heraklio, Athens, Greece}
\email{geaxatz@otenet.gr, gea.xatz@aspete.gr}
\author{I. P. STAVROULAKIS}
\address{Department of Mathematics\\
University of Ioannina\\
451 10 Ioannina, Greece}
\email{ipstav@cc.uoi.gr}

\begin{abstract}
We consider difference equations with several
non-monotone deviating arguments and nonnegative coefficients.
The deviations (delays and advances) are, generally, unbounded.
Sufficient oscillation conditions are obtained in an explicit iterative form.
Additional results in terms of $\liminf $ are obtained for bounded deviations. 
Examples illustrating the oscillation tests are presented.


\end{abstract}

\maketitle

\section{INTRODUCTION}

\footnotetext[1]{Corresponding author. E-mail maelena@math.ucalgary.ca, phone (403)-220-3956, fax 
(403)-282-5150} 

The paper deals with the difference equation with several variable retarded
arguments of the form%
\begin{equation}
\Delta x(n)+\sum_{i=1}^{m}p_{i}(n)x(\tau _{i}(n))=0\text{, \ \ \ }n\in 
\mathbb{N}
_{0}\text{,}  \tag{E$_{\text{R}}$}
\end{equation}%
where $\mathbb{N}_{0}$ is the set of nonnegative integers, and the (dual)
difference equation with several variable advanced arguments of the form%
\begin{equation}
\nabla x(n)-\sum_{i=1}^{m}p_{i}(n)x(\sigma _{i}(n))=0\text{, \ \ \ }n\in 
\mathbb{N}
\text{.}  \tag{E$_{\text{A}}$}
\end{equation}%
Equations (E$_{\text{R}}$) and (E$_{\text{A}}$) are studied under the
following assumptions: everywhere $(p_{i}(n))$, $1\leq i\leq m$, are
sequences of {\em nonnegative real numbers}, $(\tau _{i}(n))$, $1\leq i\leq m$,
are sequences of integers such that either 
\begin{equation}
\tau _{i}(n)\leq n-1,\text{ \ \ }\forall n\in 
\mathbb{N}
_{0}\text{, \ \ \ and \ \ \ }\lim\limits_{n\rightarrow \infty }\tau
_{i}(n)=\infty \text{, \ \ \ }1\leq i\leq m\text{,}  \tag{1.1}
\end{equation}%
or 
\begin{equation}
\tau _{i}(n)\leq n-1,~~\forall n\in 
\mathbb{N}
_{0}\text{~~ and ~~}\exists \text{ }M_{i}>0\text{ such that\ }n-\tau
_{i}(n)\leq M_{i}\text{,~~ }1\leq i\leq m  \tag{1.1$^{\prime }$}
\end{equation}%
and $(\sigma _{i}(n))$, $1\leq i\leq m$, are sequences of integers such that
either 
\begin{equation}
\sigma _{i}(n)\geq n+1,\text{ \ \ }\forall n\in 
\mathbb{N}
\text{, \ \ \ }1\leq i\leq m\text{,}  \tag{1.2}
\end{equation}%
or 
\begin{equation}
\sigma _{i}(n)\geq n+1,\forall n\in 
\mathbb{N}
\text{ and }\exists \text{ }\mu _{i}>0\text{ such that\ }\sigma
_{i}(n)-n\leq \mu _{i}\text{, }1\leq i\leq m\text{.}  \tag{1.2$^{\prime }$}
\end{equation}%
Here, $\Delta $ denotes the forward difference operator $\Delta
x(n)=x(n+1)-x(n)$ and $\nabla $ corresponds to the backward difference
operator $\nabla x(n)=x(n)-x(n-1)$.

Let us note that the condition $\tau _{i}(n)\leq n-1$ ($\sigma _{i}(n)\geq
n+1$) is not essential for $\limsup $-type tests; all the results of Section
2 apply to the case when $\tau _{i}(n)\leq n$ and there is a non-delay term $%
\tau _{j}(n)=n$. However, the requirement $\lim_{n\rightarrow \infty }\tau
_{i}(n)=\infty $ is significant, as Example 2.2 illustrates. A similar
remark applies to advanced difference equation (E$_{\text{A}}$).

In the special case $m=1$ equations (E$_{\text{R}}$) and (E$_{\text{A}}$)
reduce to the equations%
\begin{equation}
\Delta x(n)+p(n)x(\tau (n))=0\text{, \ \ \ }n\in 
\mathbb{N}
_{0}  \tag{1.3}
\end{equation}%
and%
\begin{equation}
\nabla x(n)-p(n)x(\sigma (n))=0\text{, \ \ \ }n\in 
\mathbb{N}
\text{,}  \tag{1.4}
\end{equation}%
respectively.

Set
\begin{equation*}
w=-\min_{\substack{ n\geq 0  \\ 1\leq i\leq m}}\tau _{i}(n)\text{.}
\end{equation*}%
Clearly, $w$ is a finite positive integer if (1.1) holds.

By a \textit{solution} of (E$_{\text{R}}$), we mean a sequence of real
numbers $(x(n))_{n\geq -w}$ which satisfies (E$_{\text{R}}$) for all $n\geq
0.$ It is clear that, for each choice of real numbers $c_{-w},$ $%
c_{-w+1},...,$ $c_{-1},$ $c_{0}$, there exists a unique solution $%
(x(n))_{n\geq -w}$ of (E$_{\text{R}}$) which satisfies the initial
conditions $x(-w)=c_{-w},$ $x(-w+1)=c_{-w+1},...,$ $x(-1)=c_{-1},$ $%
x(0)=c_{0}$.

By a solution of (E$_{\text{A}}$), we mean a sequence of real numbers $%
\left( x(n)\right) _{n\geq 0}$ which satisfies (E$_{\text{A}}$) for all $%
n\geq 1$.

A solution $(x(n))_{n\geq -w}$ (or $\left( x(n)\right) _{n\geq 0}$) of (E$_{%
\text{R}}$) (or (E$_{\text{A}}$)) is called\textit{\ oscillatory, }if the
terms $x(n)$ of the sequence are neither eventually positive nor eventually
negative. Otherwise, the solution is said to be \textit{nonoscillatory}. An
equation is \textit{oscillatory} if all its solutions oscillate.

In the last few decades, the oscillatory behavior and the existence of
positive solutions of difference equations with several deviating arguments
have been extensively studied, see, for example, papers [1-5,7-15] and
references cited therein. More precicely, in 1999, Zhang et al. [14] studied
the existence and nonexistence of positive solutions and in 1999 and 2001,
Tang et al. [10,11] considered (E$_{\text{R}}$) with constant delays,
introduced some new techniques to analyze generalized characteristic
equations and obtained several sharp oscillation conditions including
infinite sums. In 2002, Zhang et al. [15] presented some oscillation
criteria involving $\lim \sup $ for (E$_{\text{R}}$) with constant delays,
while Li et al. [7] investigated only (E$_{\text{A}}$) under the assumption
that the advances are constant and established a new oscillation condition
involving $\liminf $. In 2003, Wang [12] considered (E$_{\text{R}}$) with
constant delays and obtained some new oscillation criteria involving $%
\limsup $, while Luo et al. [8] studied a nonlinear difference equation with
several constant delays and established some sufficient oscillation
conditions. In 2005, Yan et al. [13] got oscillation criteria for (E$_{\text{%
R}}$) with variable delays. In 2006, Berezansky et al. [1] mostly
investigated non-oscillation but also obtained the following sufficient
oscillation test (Theorem 5.16): If 
\begin{equation*}
\limsup_{n\rightarrow \infty }\sum_{i=1}^{m}p_{i}(n)>0\mbox{~~and~~}%
\liminf_{n\rightarrow \infty }\sum_{i=1}^{m}\sum_{j=\tau (n)}^{n-1}p_{i}(j)>%
\frac{1}{e}\text{,}
\end{equation*}%
where $\tau (n)=\max_{1\leq i\leq m}\tau _{i}(n)$, $\forall n\geq 0$, then
all solutions of (E$_{\text{R}}$) oscillate. In the recent papers [3-5],
Chatzarakis et al. studied equations (E$_{\text{R}}$) and (E$_{\text{A}}$)
and presented some oscillation conditions involving $\limsup $ and $\liminf $.

However, most oscillation results presented in the previous papers, require
that the arguments be monotone increasing. While this condition is satisfied
by a variety of differential equations with variable delays, for difference
equations, due to the discrete nature of the arguments, if retarded
arguments are strictly increasing, then the deviations are eventually
constant. This is one of motivations to investigate difference equations
with non-monotone arguments. The challenge of this study is illustrated by
the fact that, according to [2, Theorem 3], there is no constant $A>0$ such
that the inequalities 
\begin{equation*}
\limsup_{n\rightarrow \infty }(n-\tau (n))p(n)>A\mbox{~~and~~}%
\liminf_{n\rightarrow \infty }\sum_{j=\tau (n)}^{n-1}p_{i}(j)>A
\end{equation*}%
guarantee oscillation of (1.3).

The paper is organized as follows. Section 2 contains oscillation conditions
for both retarded and advanced equations in terms of $\limsup $, where
deviations of the argument, generally, are not assumed to be bounded. In
Section 3 conditions in terms of $\liminf $ are obtained in the case when
deviations of the argument are bounded; however, a weaker result is valid
for unbounded delays. Section 4 illustrates the results of the present paper
with examples where oscillation could not have been established using
previously known results.

\section{UNBOUNDED DEVIATIONS OF ARGUMENTS AND ITERATIVE TESTS}

In 2011, Braverman and Karpuz [2] and in 2014, Stavroulakis [10],
established the following results for the special case of \ Eq. (1.3) under
the assumption that the argument $\tau (n)$ is non-monotone, (1.1) holds and 
$\xi (n)=\max_{0\leq s\leq n}\tau (s).$

\vskip0.2cm\textbf{Theorem 2.1 }([2, Theorem 5])\textbf{.} \ \textit{If}%
\begin{equation}
\mathit{\ }\limsup_{n\rightarrow \infty }\sum_{j=\xi (n)}^{n}p(j)%
\mathop{\displaystyle \prod }\limits_{i=\tau (j)}^{\xi (n)-1}\frac{1}{1-p(i)}%
>1,  \tag{2.1}
\end{equation}%
\textit{then all solutions of} (1.3) \textit{oscillate.}

\vskip0.4cm\textbf{Theorem 2.2 (}[9, Theorem 3.6]\textbf{). \ }\textit{%
Assume that} $\displaystyle a\boldsymbol{=}\liminf_{n\rightarrow \infty
}\sum_{i=\tau (n)}^{n-1}p(i)$ ~ and%
\begin{equation}
\limsup_{n\rightarrow \infty }\sum_{j=\xi (n)}^{n}p\left( j\right) %
\mathop{\displaystyle \prod }\limits_{i=\tau (j)}^{\xi (n)-1}\frac{1}{1-p(i)}%
>1-c(a),  \tag{2.2}
\end{equation}%
\textit{where}%
\begin{equation*}
c(a)=\left\{ 
\begin{array}{ll}
\frac{1}{2}\left( 1-a-\sqrt{1-2a-a^{2}}\right) , & \text{if~~~}0<a\leq 1/e,
\\ 
\frac{1}{4}\left( 2-3a-\sqrt{4-12a+a^{2}}\right) , & \text{if~~~}0<a\leq 6-4%
\sqrt{2}\text{~~and~~}p(n)\geq \frac{\alpha }{2}.%
\end{array}%
\right.
\end{equation*}%
\textit{Then all solutions of equation }(1.3)\textit{\ oscillate.}

\vskip0.2cmIn this section, we study oscillation of (E$_{\text{R}}$) and (E$%
_{\text{A}}$). We establish new sufficient oscillation conditions involving $%
\lim \sup $ under the assumption that the arguments are non-monotone and
(1.1) holds. Even for equation (1.3) with one delay, our results improve
(2.1) and (2.2). The method we use is based on the iterative application of
the Gronwall inequality.

\subsection{Retarded difference equations}

Let%
\begin{equation}
\varphi _{i}(n)=\max_{0\leq s\leq n}\tau _{i}(s),\text{ \ \ }n\geq 0 
\tag{2.3}
\end{equation}%
and%
\begin{equation}
\varphi (n)=\max_{1\leq i\leq m}\varphi _{i}(n),\text{ \ \ }n\geq 0\text{.} 
\tag{2.4}
\end{equation}%
Clearly, the sequences of integers $\varphi (n)$, $\varphi _{i}(n)$, $1\leq
i\leq m$ are non-decreasing and $\varphi (n)\leq n-1$, $\varphi _{i}(n)\leq
n-1$, $1\leq i\leq m$ for all $n\geq 0$.

The following simple result is cited to explain why we can consider only the
case 
\begin{equation}
\sum_{i=1}^{m}p_{i}(n)<1,\text{\ \ \ } ~~~~~\forall n\geq 0\text{.} 
\tag{2.5}
\end{equation}

\vskip0.2cm\textbf{Theorem 2.3. } \textit{Assume that there exists a
subsequence }$\theta (n)$, $n\in 
\mathbb{N}
$ \textit{of positive integers} \textit{such that}%
\begin{equation}
\sum_{i=1}^{m}p_{i}(\theta (n))\geq 1,\text{\ \ \ }~~~~~\forall n\in {%
\mathbb{N}}\text{.}  \tag{2.6}
\end{equation}%
\textit{Then all solutions of} (E$_{\text{R}}$) \textit{oscillate}.

\begin{proof}
Assume, for the sake of contradiction, that $\left( x(n)\right) _{n\geq -w}$
is a nonoscillatory solution of (E$_{\text{R}}$). Then it is either
eventually positive or eventually negative. As $\left( -x(n)\right) _{n\geq
-w}$ is also a solution of (E$_{\text{R}}$), we may restrict ourselves only
to the case where $x(n)>0$ for all large $n.$ Let $n_{1}\geq -w$ be an
integer such that $x(n)>0$ for all $n\geq n_{1}$. Then, there exists $%
n_{2}\geq n_{1}$ such that $x(\tau _{i}(n))>0$,\ $\forall n\geq n_{2}$, $%
1\leq i\leq m$. In view of this, Eq.(E$_{\text{R}}$) becomes%
\begin{equation*}
\Delta x(n)=-\sum_{i=1}^{m}p_{i}(n)x(\tau _{i}(n))\leq 0\text{, \ \ \ \ }%
\forall \theta(n)\geq n_{2}\text{,}
\end{equation*}%
which means that the sequence $(x(n))$ is eventually decreasing.

Taking into account the fact that (2.6) holds, equation (E$_{\text{R}}$)
gives%
\begin{eqnarray*}
x(\theta (n)+1) &=&x(\theta (n))-\sum_{i=1}^{m}p_{i}(\theta (n))x(\tau
_{i}(\theta (n)))\leq x(\theta (n))-x(\theta (n))\sum_{i=1}^{m}p_{i}(\theta
(n)) \\
&=&x(\theta (n))\left( 1-\sum_{i=1}^{m}p_{i}(\theta (n))\right) \leq 0\text{%
, \ \ \ \text{for all}}\mathit{\ }\theta (n)\geq n_{2},
\end{eqnarray*}%
where $\theta (n)\rightarrow \infty $ as $n\rightarrow \infty $, which
contradicts the assumption that $x(n)>0$ for all $n\geq n_{1}$. 
\end{proof}

In a series of further computations, we try to extend the iterative result,
established by Koplatadze and Kvinikadze in [6] for a differential equation
with a single delay, to difference equations with several retarded
arguments. The following lemma provides an estimation of a positive solution
rate of decay, which is a useful tool for obtaining oscillation conditions.

\vskip0.2cm\textbf{Lemma 2.1. } \textit{Assume that }(2.5) \textit{holds and 
}$x(n)$ \textit{is a positive solution of }(E$_{\text{R}}$).\textit{\ Set}%
\begin{equation}
a_{1}(n,k):=\prod_{i=k}^{n-1}\left[ 1-\sum_{\ell =1}^{m}p_{\ell }(i)\right] 
\tag{2.7}
\end{equation}%
\textit{and}%
\begin{equation}
a_{r+1}(n,k):=\prod_{i=k}^{n-1}\left[ 1-\sum_{\ell =1}^{m}p_{\ell
}(i)a_{r}^{-1}(i,\tau _{\ell }(i))\right] ,\text{ \ \ \ \ }r\in \mathbb{N}. 
\tag{2.8}
\end{equation}%
\textit{Then}%
\begin{equation}
x(n)\leq a_{r}(n,k)x(k),\text{ \ \ \ \ }r\in \mathbb{N}.  \tag{2.9}
\end{equation}

\begin{proof}
Since $x(n)$ is a positive solution of (E$_{\text{R}}$), 
\begin{equation*}
\Delta x(n)=-\sum_{i=1}^{m}p_{i}(n)x(\tau _{i}(n))\leq 0\text{, \ \ \ \ }
\forall n\in {\mathbb{N}}_{0}\text{,}
\end{equation*}%
which means that the sequence $(x(n))$ is decreasing.

From (E$_{\text{R}}$) and the decreasing character of $(x(n))$ , we have%
\begin{equation*}
\Delta x(n)+x(n)\sum_{i=1}^{m}p_{i}(n)\leq 0\text{, \ \ \ \ } \forall n\in 
\mathbb{N}_{0}\text{.}
\end{equation*}%
Applying the discrete Gronwall inequality, we obtain%
\begin{equation*}
x(n)\leq x(k)\mathop{\displaystyle \prod }\limits_{i=k}^{n-1}\left[
1-\sum_{\ell =1}^{m}p_{\ell }(i)\right] =a_{1}(n,k)x(k),\text{ \ \ \ \ }
\forall n\geq k\geq 0\text{,}
\end{equation*}%
which justifies (2.9) for $r=1$.

Assume that (2.9) holds for some $r>1$. Substituting%
\begin{equation*}
x(\tau _{\ell }(n))\geq x(n)a_{r}^{-1}(n,\tau _{\ell }(n))
\end{equation*}%
into (E$_{\text{R}}$) leads to the inequality%
\begin{equation*}
\Delta x(n)+\sum_{\ell =1}^{m}p_{\ell }(n)x(n)a_{r}^{-1}(n,\tau _{\ell
}(n))\leq 0\text{, \ \ \ for sufficiently large }n\text{.}
\end{equation*}%
Again, applying the discrete Gronwall inequality, we obtain%
\begin{equation*}
x(n)\leq x(k)\mathop{\displaystyle \prod }\limits_{i=k}^{n-1}\left[
1-\sum_{\ell =1}^{m}p_{\ell }(i)a_{r}^{-1}(i,\tau _{\ell }(i))\right]
=a_{r+1}(n,k)x(k),
\end{equation*}%
due to the definition of $a_{r+1}$ in (2.8), which concludes the induction
step. The proof of the lemma is complete.
\end{proof}

It is interesting to see how the estimate developed in Lemma 2.1 works in
the case of autonomous equations.

\vskip0.2cm\textbf{Example 2.1. \ } For the equation 
\begin{equation*}
\Delta x(n)+\frac{1}{4}x(n-1)=0,
\end{equation*}%
inequality (2.9) for $k=n-1$ has the form 
\begin{equation*}
x(n)\leq a_{r}(n,n-1)x(n-1).
\end{equation*}%
Here $a_{1}(n,n-1)=\frac{3}{4}$, $a_{2}(n,n-1)=\frac{2}{3}$, $a_{3}(n,n-1)=%
\frac{5}{8}$, $a_{4}(n,n-1)=\frac{3}{5}$, $a_{5}(n,n-1)=\frac{7}{12}$, $%
\dots ,$ $a_{2k-1}=\frac{2k+1}{4k}$, $a_{2k}=\frac{k+1}{2k+1}$, thus the
best possible estimate of the rate of decrease of a nonoscillatory solution
is $x(n+1)\leq \frac{1}{2}x(n)$, or $x(n)\leq x(0) (\frac{1}{2})^{n}$, which
is in compliance with the fact that $\lambda =\frac{1}{2}$ is the (double)
root of the characteristic equation $\lambda ^{2}-\lambda +\frac{1}{4}=0$.

\vskip0.2cm\textbf{Theorem 2.4. } \textit{Assume that }$(p_{i}(n))$, $1\leq
i\leq m$, \textit{are sequences of nonnegative real numbers, }(1.1) \textit{%
and }(2.5) \textit{hold, }$\varphi (n)$ \textit{is defined by }(2.4) and
$a_{r}(n,k)$ are denoted in (2.7),(2.8).
\textit{If}
\begin{equation}
\limsup_{n\rightarrow \infty }\sum_{j=\varphi
(n)}^{n}\sum_{i=1}^{m}p_{i}(j)a_{r}^{-1}(\varphi (n),\tau _{i}(j))>1, 
\tag{2.10}
\end{equation}%
\textit{then all solutions of} (E$_{\text{R}}$) \textit{oscillate}.

\begin{proof}
Assume, for the sake of contradiction, that $(x(n))_{n\geq -w}$ is a
nonoscillatory solution of (E$_{\text{R}}$). Then it is either eventually
positive or eventually negative. Similarly to the proof of Theorem 2.3, we
may restrict ourselves only to the case $x(n)>0$ for all large $n.$ Let $%
n_{1}\geq -w$ be an integer such that $x(n)>0$ for all $n\geq n_{1}$. By
(1.1), there exists $n_{2}\geq n_{1} $ such that $x(\tau _{i}(n))>0$,\ $%
\forall n\geq n_{2}$, $1\leq i\leq m $. In view of this, Eq.(E$_{\text{R}}$)
becomes%
\begin{equation*}
\Delta x(n)=-\sum_{i=1}^{m}p_{i}(n)x(\tau _{i}(n))\leq 0\text{, \ \ \ \ }%
\forall n\geq n_{2}\text{,}
\end{equation*}%
which means that the sequence $(x(n))$ is eventually decreasing.

Summing up (E$_{\text{R}}$) from $\varphi (n)$ to $n$, and using the fact
that the function $x$ is non-increasing, while the function $\varphi $ (as
defined by (2.4)) is non-decreasing, and taking into account that%
\begin{equation*}
\tau _{i}(j)\leq \varphi (n)\text{ \ \ and \ \ }x(\tau _{i}(j))\geq
x(\varphi (n))a_{r}^{-1}(\varphi (n),\tau _{i}(j))\text{,}
\end{equation*}%
we obtain, for sufficiently large $n$,%
\begin{eqnarray*}
x(\varphi (n)) &=&x(n+1)+\sum_{j=\varphi (n)}^{n}\sum_{i=1}^{m}p_{i}\left(
j\right) x(\tau _{i}(j)) \\
&\geq &\sum_{j=\varphi (n)}^{n}\sum_{i=1}^{m}p_{i}\left( j\right) x(\tau
_{i}(j))\geq x(\varphi (n))\sum_{j=\varphi
(n)}^{n}\sum_{i=1}^{m}p_{i}(j)a_{r}^{-1}(\varphi (n),\tau _{i}(j))\text{.}
\end{eqnarray*}%
Consequently,%
\begin{equation*}
x(\varphi (n))\left[ 1-\sum_{j=\varphi
(n)}^{n}\sum_{i=1}^{m}p_{i}(j)a_{r}^{-1}(\varphi (n),\tau _{i}(j))\right]
\leq 0\text{,}
\end{equation*}%
which gives%
\begin{equation*}
\limsup_{n\rightarrow \infty }\sum_{j=\varphi
(n)}^{n}\sum_{i=1}^{m}p_{i}(j)a_{r}^{-1}(\varphi (n),\tau _{i}(j))\leq 1%
\text{.}
\end{equation*}%
This contradicts (2.10) and the proof of the theorem is complete.
\end{proof}

\vskip0.1cmTo establish the next theorem we need the following lemma.

\vskip0.2cm\textbf{Lemma 2.2. \ }(cf. [3, Lemma 2.1])\textbf{. \ }\textit{%
Assume that }(1.1) \textit{holds, }$\left( x(n)\right) $ \textit{is an
eventually positive solution of }(E$_{\text{R}}$)\textit{, and}%
\begin{equation}
\alpha =\min \left\{ \alpha _{i}:1\leq i\leq m\right\} \text{,}  \tag{2.11}
\end{equation}%
\textit{where}%
\begin{equation}
\alpha _{i}=\liminf_{n\rightarrow \infty }\sum_{j=\varphi
_{i}(n)}^{n-1}p_{i}\left( j\right) \text{.}  \tag{2.12}
\end{equation}

\textit{If }$0<\alpha \leq -1+\sqrt{2}$,\textit{\ then}%
\begin{equation}
\liminf_{n\rightarrow \infty }\frac{x(n+1)}{x(\varphi (n))}\geq \frac{%
1-\alpha -\sqrt{1-2\alpha -\alpha ^{2}}}{2}\text{.}  \tag{2.13}
\end{equation}

\vskip0.1cm\textbf{Theorem 2.5. \ }\textit{Assume that }$(p_{i}(n))$, $1\leq
i\leq m$, \textit{are sequences of nonnegative real numbers, }(1.1) \textit{%
and }(2.5) \textit{hold, and }$\varphi (n)$ \textit{is defined by }(2.4).%
\textit{\ If }$0<\alpha \leq 1/e$\textit{,} \textit{where }$\alpha $ \textit{%
is denoted by (2.11), and}%
\begin{equation}
\limsup_{n\rightarrow \infty }\sum_{j=\varphi
(n)}^{n}\sum_{i=1}^{m}p_{i}(j)a_{r}^{-1}(\varphi (n),\tau _{i}(j))>1-\frac{%
1-\alpha -\sqrt{1-2\alpha -\alpha ^{2}}}{2}\text{,}  \tag{2.14}
\end{equation}%
\textit{where $a_{r}(n,k)$ is as in }(2.7), (2.8),\textit{\ then all
solutions of} (E$_{\text{R}}$) \textit{oscillate}.

\vskip0.1cm\textbf{Proof.\ \ }Assume, for the sake of contradiction, that $%
(x(n))_{n\geq -w}$ is a nonoscillatory solution of (E$_{\text{R}}$). Then,
as in the proof of Theorem 2.4, for sufficiently large $n$, we obtain%
\begin{eqnarray*}
x(\varphi (n)) &=&x(n+1)+\sum_{j=\varphi (n)}^{n}\sum_{i=1}^{m}p_{i}\left(
j\right) x(\tau _{i}(j)) \\
&\geq &x(n+1)+x(\varphi (n))\sum_{j=\varphi
(n)}^{n}\sum_{i=1}^{m}p_{i}(j)a_{r}^{-1}(\varphi (n),\tau _{i}(j))\text{.}
\end{eqnarray*}%
that is,%
\begin{equation*}
\sum_{j=\varphi (n)}^{n}\sum_{i=1}^{m}p_{i}(j)a_{r}^{-1}(\varphi (n),\tau
_{i}(j))\leq 1-\frac{x(n+1)}{x(\varphi (n))}\text{, }
\end{equation*}%
which gives%
\begin{equation*}
\limsup_{n\rightarrow \infty }\sum_{j=\varphi
(n)}^{n}\sum_{i=1}^{m}p_{i}(j)a_{r}^{-1}(\varphi (n),\tau _{i}(j))\leq
1-\liminf_{n\rightarrow \infty }\frac{x(n+1)}{x(\varphi (n))}\text{.}
\end{equation*}%
Assume that $0<\alpha \leq 1/e$ (clearly, $\alpha <-1+\sqrt{2})$ and, by
Lemma 2.2, inequality (2.13) holds, and so the last inequality leads to%
\begin{equation*}
\limsup_{n\rightarrow \infty }\sum_{j=\varphi
(n)}^{n}\sum_{i=1}^{m}p_{i}(j)a_{r}^{-1}(\varphi (n),\tau _{i}(j))\leq 1-%
\frac{1-\alpha -\sqrt{1-2\alpha -\alpha ^{2}}}{2}\text{,}
\end{equation*}%
which contradicts condition (2.14).

The proof of the theorem is complete.

\vskip0.2cm\textbf{Remark 2.1.} \ Observe that conditions (2.1) and (2.2)
are special cases of (2.10) and (2.14) respectively, when $r=1$.

\vskip0.2cm\textbf{Remark 2.2.} \ The following example illustrates the
significance of the condition $\lim\limits_{n\rightarrow \infty }\tau
_{i}(n)=\infty $, $1\leq i\leq m$, in Theorems 2.4 and 2.5.

\vskip0.2cm\textbf{Example 2.2.} \ Consider the retarded difference equation
(1.3) with 
\begin{equation*}
p(n)=\left\{ 
\begin{array}{ll}
\displaystyle1/2, & \text{\ \ if }n=3k, \\ 
\displaystyle3/10, & \text{\ \ if }n\neq 3k,%
\end{array}%
\right. \text{ \ \ }\tau (n)=\left\{ 
\begin{array}{ll}
-1, & \text{\ \ if }n=3k, \\ 
n-1, & \text{\ \ if }n\neq 3k,%
\end{array}%
\right. \text{\ \ \ }k\in \mathbb{N}_{0}.
\end{equation*}%
Obviously (2.5) is satisfied. Also, by (2.4),%
\begin{equation*}
\varphi (n)=\max_{0\leq s\leq n}\tau (s)=\left\{ 
\begin{array}{ll}
n-2, & \text{\ \ if }n=3k, \\ 
n-1, & \text{\ \ if }n\neq 3k,%
\end{array}%
\right. k\in \mathbb{N}_{0}.
\end{equation*}%
If $n=3k$, then%
\begin{equation*}
\sum_{j=\varphi (n)}^{n}p(j)=\sum_{j=n-2}^{n}p(j)=\frac{3}{10}+\frac{3}{10}+%
\frac{1}{2}=\frac{11}{10}.
\end{equation*}%
Therefore $\displaystyle\limsup_{n\rightarrow \infty }\sum_{j=\varphi
(n)}^{n}p_{i}(j)\geq \frac{11}{10}>1,$ which means that (2.10) is satisfied
for any $r$.

Also, if $n=3k$, then%
\begin{equation*}
\sum_{j=\varphi (n)}^{n-1}p_{i}\left( j\right)
=\sum_{j=n-2}^{n-1}p_{i}\left( j\right) =\frac{3}{10}+\frac{3}{10}=\frac{3}{5%
}
\end{equation*}%
and, if $n\neq 3k$, then%
\begin{equation*}
\sum_{j=\varphi (n)}^{n-1}p(j)=\sum_{j=n-1}^{n-1}p(j)=\left\{ 
\begin{array}{ll}
1/2, & \text{\ \ if }n=3k+1, \\ 
3/10, & \text{\ \ if }n=3k+2,%
\end{array}%
\right. k\in \mathbb{N}_{0}.
\end{equation*}%
Therefore%
\begin{equation*}
\alpha =\liminf_{n\rightarrow \infty }\sum_{j=\varphi (n)}^{n-1}p_{i}\left(
j\right) =\min \left\{ \frac{3}{5},\frac{1}{2},\frac{3}{10}\right\} =\frac{3%
}{10}<\frac{1}{e}
\end{equation*}%
and%
\begin{equation*}
\sum_{j=\varphi (n)}^{n}p_{i}(j)=\frac{11}{10}>1-\frac{1-\alpha -\sqrt{%
1-2\alpha -\alpha ^{2}}}{2}\simeq 0.928388218\text{,}
\end{equation*}%
which means that (2.14) is satisfied for any $r$. Observe, however, that
equation (1.3) has a nonoscillatory solution 
\begin{equation*}
\begin{array}{ll}
x(-1)=-\frac{12}{7}, & x(0)=\displaystyle1,~x(1)=\frac{13}{7},~x(2)=\frac{109%
}{70}, \\ 
&  \\ 
& \displaystyle x(3)=1,~x(4)=\frac{13}{7},~x(5)=\frac{109}{70}%
,~x(6)=1,~~\dots ~%
\end{array}%
\end{equation*}%
which illustrates the significance of the condition $\lim\limits_{n%
\rightarrow \infty }\tau (n)=\infty $ in Theorems 2.4 and 2.5.

\subsection{Advanced difference equations}

Similar oscillation theorems for the (dual) advanced difference equation (E$%
_{\text{A}}$) can be derived easily. The proofs of these theorems are
omitted, since they follow a similar procedure as in Subsection 2.1.

Denote 
\begin{equation}
\rho _{i}(n)=\min_{s\geq n}\sigma _{i}(s),\text{ \ \ }n\geq 0 
\tag{2.15}
\end{equation}%
and%
\begin{equation}
\rho (n)=\min_{1\leq i\leq m}\rho _{i}(n),\text{ \ \ }n\geq 0\text{.} 
\tag{2.16}
\end{equation}
Clearly, the sequences of integers $\rho (n)$, $\rho _{i}(n)$, $1\leq i\leq
m $ are non-decreasing and $\rho (n)\geq n+1$, $\rho _{i}(n)\geq n+1$, $%
1\leq i\leq m$ for all $n\geq 0$.

\vskip0.2cm\textbf{Theorem 2.3}$^{\prime }$\textbf{. } \textit{Assume that
there exists a subsequence }$\theta (n),~n\in 
\mathbb{N}
$ \textit{of positive integers such that}%
\begin{equation*}
\sum_{i=1}^{m}p_{i}(\theta (n))\geq 1\mathit{\ \ \ \ \ }\forall n\in {%
\mathbb{N}}\text{.}
\end{equation*}%
\textit{Then all solutions of} (E$_{\text{A}}$) \textit{oscillate}.

\vskip0.2cm\textbf{Theorem 2.4}$^{\prime }$. \ \textit{Assume that }$%
(p_{i}(n))$, $1\leq i\leq m$, \textit{are sequences of nonnegative real
numbers, }(1.2) \textit{and }(2.5) \textit{hold,} $\rho (n)$ \textit{is
defined by }(2.16) and%
\begin{equation}
b_{1}(n,k):=\prod_{i=n+1}^{k}\left[ 1-\sum_{\ell =1}^{m}p_{\ell }(i)\right] 
\tag{2.17}
\end{equation}%
and%
\begin{equation}
b_{r+1}(n,k):=\prod_{i=n+1}^{k}\left[ 1-\sum_{\ell =1}^{m}p_{\ell
}(i)b_{r}^{-1}(i,\sigma _{\ell }(i))\right] ,\text{ \ \ \ \ }r\in \mathbb{N}.
\tag{2.18}
\end{equation}%
\textit{If}%
\begin{equation}
\limsup_{n\rightarrow \infty }\sum_{j=n}^{\rho
(n)}\sum_{i=1}^{m}p_{i}(j)b_{r}^{-1}(\rho (n),\sigma _{i}(j))>1\text{,} 
\tag{2.19}
\end{equation}%
\textit{then all solutions of} (E$_{\text{A}}$) \textit{oscillate}.

\vskip0.2cm\textbf{Lemma 2.2}$^{\prime }$ (cf. [3, Lemma 2.1]). \ \textit{%
Assume that }(1.2) \textit{holds, }$\left( x(n)\right) $ \textit{is a n
eventually positive solution of }(E$_{\text{R}}$)\textit{, and}%
\begin{equation}
\alpha =\min \left\{ \alpha _{i}:1\leq i\leq m\right\} \text{,}  \tag{2.20}
\end{equation}%
\textit{where}%
\begin{equation}
\alpha _{i}=\liminf_{n\rightarrow \infty }\sum_{j=n+1}^{\rho
_{i}(n)}p_{i}\left( j\right) \text{.}  \tag{2.21}
\end{equation}

\textit{If }$0<\alpha \leq -1+\sqrt{2}$,\textit{\ then}%
\begin{equation}
\liminf_{n\rightarrow \infty }\frac{x(n-1)}{x(\rho (n))}\geq \frac{1-\alpha -%
\sqrt{1-2\alpha -\alpha ^{2}}}{2}\text{.}  \tag{2.22}
\end{equation}

\vskip0.2cm\textbf{Theorem \ 2.5}$^{\prime }$\textbf{. }\ \textit{Assume
that }$(p_{i}(n))$, $1\leq i\leq m$ \textit{are sequences of nonnegative
real numbers, }(1.2) \textit{and }(2.5) \textit{hold, }$\rho (n)$ \textit{is
defined by }(2.16)\textit{, and }$a$ \textit{is denoted by }(2.20)\textit{.
If }$0<\alpha \leq 1/e$\textit{,} \textit{and}%
\begin{equation}
\limsup_{n\rightarrow \infty }\sum_{j=n}^{\rho
(n)}\sum_{i=1}^{m}p_{i}(j)b_{r}^{-1}(\rho (n),\sigma _{i}(j))>1-\frac{%
1-\alpha -\sqrt{1-2\alpha -\alpha ^{2}}}{2}\text{,}  \tag{2.23}
\end{equation}%
\textit{where $b_{r}(n,k)$ is defined in }(2.17), (2.18),\textit{\ then all
solutions of} (E$_{\text{A}}$) \textit{oscillate}.

\section{BOUNDED DEVIATIONS OF ARGUMENTS}

In this section, new sufficient oscillation conditions, involving $\lim \inf 
$, under the assumption that all the delays (advances) are bounded and (1.1$%
^{\prime }$) $\left( \text{or (1.2}^{\prime }\text{)}\right) $ holds, are
established for equations (E$_{\text{R}}$) and (E$_{\text{A}}$). These
conditions improve Theorem 5.16 in [1] and extend the following oscillation
results by Chatzarakis et al. [5] to the case of non-monotone arguments.

\vskip0.2cm\textbf{Theorem 3.1 (}[5, Theorem 2.2]\textbf{). \ }\textit{%
Assume that the sequences }$\left( \tau _{i}(n)\right) $\textit{, }$1\leq
i\leq m$ \textit{are increasing, }(1.1)\textit{\ holds,}%
\begin{equation}
\limsup_{n\rightarrow \infty }\sum_{i=1}^{m}p_{i}(n)>0\text{\ \ \ \ and \ \
\ \ }\liminf_{n\rightarrow \infty }\sum_{i=1}^{m}\sum_{j=\tau
_{i}(n)}^{n-1}p_{i}(j)>\frac{1}{e}\text{\thinspace .}  \tag{3.1}
\end{equation}%
\textit{Then all solutions of Eq. }(E$_{\text{R}}$)\textit{\ oscillate.}

\vskip0.2cm\textbf{Theorem 3.2 (}[5, Theorem 3.2]\textbf{). \ }\textit{%
Assume that the sequences }$\left( \sigma _{i}(n)\right) $\textit{,\ }$1\leq
i\leq m$ \textit{are increasing, }(1.2)\textit{\ holds,}%
\begin{equation}
\limsup_{n\rightarrow \infty }\sum_{i=1}^{m}p_{i}(n)>0\text{\ \ \ \ and \ \
\ \ }\liminf_{n\rightarrow \infty }\sum_{i=1}^{m}\sum_{j=n+1}^{\sigma
_{i}(n)}p_{i}(j)>\frac{1}{e}\text{\thinspace .}  \tag{3.2}
\end{equation}%
\textit{Then all solutions of Eq. }(E$_{\text{A}}$)\textit{\ oscillate.}

\subsection{Retarded difference equations}

We present a new sufficient oscillation condition for (E$_{\text{R}}$),
under the assumption that all the delays are bounded and (1.1$^{\prime }$)
holds.


\vskip0.2cm\textbf{Theorem 3.3. } \ \textit{Assume that }$(p_{i}(n))$, $%
1\leq i\leq m$, \textit{are sequences of nonnegative real numbers,} (1.1$%
^{\prime }$) \textit{holds and for each} $k=1, \dots, m$,
\begin{equation}
\liminf_{n\rightarrow \infty }  
\sum_{i=1}^{m} 
\sum_{j=\tau_{k}(n)}^{n-1} 
p_{i}(j)>\left( \frac{M_k}{M_k+1}\right) ^{M_k+1}\text{,}  \tag{3.3}
\end{equation}%
\textit{then all solutions of} (E$_{\text{R}}$) \textit{oscillate.}

\begin{proof}
Assume, for the sake of contradiction, that $(x(n))_{n\geq -w}$ is a
nonoscillatory solution of (E$_{\text{R}}$). Without loss of generality, we
can assume that $x(n)>0$ for all $n\geq n_{1}$. Then, there exists $n_{2}\geq n_{1}$
such that $x(\tau _{i}(n))>0$,\ $\forall n\geq n_{2}$, $1\leq i\leq m$. In
view of this, Eq.(E$_{\text{R}}$) becomes%
\begin{equation*}
\Delta x(n)=-\sum_{i=1}^{m}p_{i}(n)x(\tau _{i}(n))\leq 0\text{, \ \ \ \ }%
\forall n\geq n_{2}\text{,}
\end{equation*}%
which means that the sequence $(x(n))_{n\geq n_{2}}$ is non-increasing. The
sequences%
\begin{equation*}
b_{k}(n)=\left( \frac{n-\tau _{k}(n)}{n-\tau _{k}(n)+1}\right) ^{n-\tau
_{k}(n)+1}
\end{equation*}%
satisfy the inequality%
\begin{equation}
\frac{1}{4}\leq b_{k}(n)\leq \left( \frac{M_k}{M_k+1}\right)^{M_k+1},\text{ \ \ \ 
}n\geq 1,\text{ \ \ \ \ }1\leq k\leq m.  \tag{3.4}
\end{equation}%
Due to (3.3), for each $k=1, \dots,m$, we can choose $n_{0}(k) \geq n_{2}$ and 
$\varepsilon_{k}>0$ such that for $n \geq n_{0}(k)$,
$$
\sum_{i=1}^{m}\sum_{j=\tau_k(n)}^{n-1} p_{i}(j)>\left( \frac{M_k}{M_k+1}%
\right) ^{M_k+1}+\varepsilon_{k}
$$
and
$$
d_k :=\left( \frac{M_k}{M_k+1}\right) ^{-M_k-1}
\left[ \left( \frac{M_k}{M_k+1}\right)^{M_k+1}+\varepsilon _{k}\right] > 1.
$$
Denoting
$$
\varepsilon_{0} := \min_{1\leq k \leq m} \varepsilon_{k}>0, 
\quad d := \min_{1\leq k \leq m} d_k>1,
$$
we obtain
\begin{equation*}
\sum_{i=1}^{m}\sum_{j=\tau_{k}(n)}^{n-1}p_{i}(j)>\left( \frac{M_k}{M_k+1}%
\right) ^{M_k+1}+\varepsilon _{0}
\end{equation*}
and
\begin{equation*}
\left( \frac{M_k}{M_k+1}\right) ^{-M_k-1}\left[ \left( \frac{M_k}{M_k+1}\right)
^{M_k+1}+\varepsilon_{0}\right] \geq d >1,
\end{equation*}%
which together with (3.3) immediately implies for any $k=1,\dots ,m$
\begin{equation*}
\sum_{i=1}^{m}\sum_{j=\tau _{k}(n)}^{n-1}\frac{p_{i}(j)}{b_{k}(n)}\geq
\left( \frac{M_k}{M_k+1}\right) ^{-M_k-1}\sum_{i=1}^{m}\sum_{j=\tau_{k}(n)}^{n-1} p_{i}(j) > d>1.
\end{equation*}
Dividing (E$_{\text{R}}$) by $x(n)$ we have 
\begin{equation*}
\frac{x(n+1)}{x(n)}=1-\sum_{i=1}^{m}p_{i}(n)\frac{x(\tau _{i}(n))}{x(n)},%
\text{ \ \ \ \ }n\geq n_{0}.
\end{equation*}%
Multiplying and taking into account that $x(\tau _{i}(n))/x(n)\geq 1$, we
obtain the estimate%
\begin{equation*}
\frac{x(n)}{x(\tau _{k}(n))}=\prod_{j=\tau _{k}(n)}^{n-1}\frac{x(j+1)}{x(j)}%
\leq \prod_{j=\tau _{k}(n)}^{n-1}\left( 1-\sum_{i=1}^{m}p_{i}(j)\right) .
\end{equation*}%
Using this estimate and the relation between the arithmetic and the
geometric means, we have 
\begin{equation}
\frac{x(n)}{x(\tau _{k}(n))}\leq \left[ 1-\frac{1}{n-\tau _{k}(n)}%
\sum_{i=1}^{m}\sum_{j=\tau _{k}(n)}^{n-1}p_{i}(j)\right] ^{n-\tau _{k}(n)}. 
\tag{3.5}
\end{equation}%
Observe that the function $f:\left( 0,1\right) \rightarrow 
\mathbb{R}
$ defined as%
\begin{equation*}
f(y):=y(1-y)^{\rho },\text{ \ \ \ \ }\rho \in {\mathbb{N}},
\end{equation*}%
attains its maximum at $y=\frac{1}{1+\rho }$, which equals $\displaystyle %
f_{\max }=\frac{\rho ^{\rho }}{(1+\rho )^{1+\rho }}\text{\thinspace .}$ In
the inequality 
\begin{equation*}
y(1-y)^{\rho }\leq \frac{\rho ^{\rho }}{(1+\rho )^{1+\rho }},\text{ \ \ \ \ }%
y\in (0,1),\text{ \ \ \ \ }\rho \in {\mathbb{N}}\text{,}
\end{equation*}%
assuming $\rho =n-\tau _{k}(n)$, $y=x/\rho$, where $\displaystyle 
x=\sum_{i=1}^{m}\sum_{j=\tau_{k}(n)}^{n-1}p_{i}(j)$, we obtain from (3.5) 
\begin{equation*}
\frac{x(\tau_{k}(n))}{x(n)}\geq \sum_{i=1}^{m}\sum_{j=\tau
_{k}(n)}^{n-1}p_{i}(j)\left( \frac{n-\tau_{k}(n)+1}{n-\tau _{k}(n)}\right)
^{n-\tau _{k}(n)+1}=\sum_{i=1}^{m}\sum_{j=\tau_{k}(n)}^{n-1}\frac{p_{i}(j)}{b_{k}(n)}>d
\end{equation*} 
for any $n\geq n_{0}(k)$. Denote $n_0=\max\{ n_0(1), \dots, n_0(m) \}$. 
If we continue this procedure assuming $n\geq n_{0}+M$,
where $M=\max_{1\leq i\leq m}M_{i}$,
and using the properties of the geometric and the algebraic mean, we have 
\begin{equation*}
\frac{x(n)}{x(\tau _{k}(n))}=\prod_{j=\tau_{k}(n)}^{n-1}\frac{x(j+1)}{x(j)}
=\prod_{j=\tau _{k}(n)}^{n-1}\left( 1-\sum_{i=1}^{m}p_{i}(j)\frac{x(\tau_{i}(j))}{x(j)}\right)
\end{equation*}
\begin{equation*}
\leq \prod_{j=\tau _{k}(n)}^{n-1}\left( 1-d\sum_{i=1}^{m}p_{i}(j)\right)
\leq \left[ 1-\frac{d}{n-\tau 
_{k}(n)}\sum_{i=1}^{m}\sum_{j=\tau_{k}(n)}^{n-1}p_{i}(j)\right]^{n-\tau _{k}(n)}.
\end{equation*}%
Applying the same argument, we obtain 
\begin{equation*}
\frac{x(\tau _{k}(n))}{x(n)}\geq d\sum_{i=1}^{m}\sum_{j=\tau_{k}(n)}^{n-1}
\frac{p_{i}(j)}{b_{k}(n)}>d^{2},\text{ \ }~n\geq n_{0}+2M,~~k=1,\dots ,m,
\end{equation*}
\begin{equation*}
\frac{x(\tau_{k}(n))}{x(n)}>d^{r},\text{ \ \ \ }n\geq n_{0}+ rM,\text{ \ \ \ 
}k=1,\dots ,m.
\end{equation*}%
Due to (3.3), we observe that for any $k$
$$
\limsup_{n\rightarrow \infty }\sum_{i=1}^{m}p_{i}(n)\geq \frac{1}{M_k}
\left( \frac{M_k}{M_k+1} \right)^{M_k+1}. 
$$
Since the function $\displaystyle h(x) := \frac{1}{x}  
\frac{x^{x+1}}{(x+1)^{x+1}}=\frac{x^x}{(x+1)^{x+1}}$, $x \geq 1$, is decreasing and $M_k \leq 
M$,  we conclude that
\begin{equation}
\limsup_{n\rightarrow \infty }\sum_{i=1}^{m}p_{i}(n)\geq c:=\frac{1}{M}
\left( \frac{M}{M+1}\right) ^{M+1}  \tag{3.6}
\end{equation}
and choose a subsequence $(\theta (n))$ of ${\mathbb{N}}$ such that 
\begin{equation*}
\sum_{i=1}^{m}p_{i}(\theta (n))\geq c>0.
\end{equation*}%
Since 
\begin{equation*}
0<\frac{x(n+1)}{x(n)}=1-\sum_{i=1}^{m}p_{i}(n)\frac{x(\tau _{i}(n))}{x(n)}~,
\end{equation*}
we have 
\begin{equation*}
\sum_{i=1}^{m}p_{i}(n)\frac{x(\tau _{i}(n))}{x(n)}<1.
\end{equation*}%
In particular, 
\begin{equation*}
\min_{1\leq k\leq m}\frac{x(\tau _{k}(\theta (n)))}{x(\theta (n))}%
\sum_{i=1}^{m}p_{i}(\theta (n))<1.
\end{equation*}%
Choosing $r \in {\mathbb{N}}$ such that $d^{r}>\frac{1}{c}$, where $c$ was
defined in (3.6), $\theta(n)\geq n_{0}+ rM$ and noticing that 
\begin{equation*}
d^{r}<\min_{1\leq k\leq m}\frac{x(\tau_{k}(\theta (n)))}{x(\theta (n))}\leq
\left( \sum_{i=1}^{m}p_{i}(\theta (n))\right)^{-1}\leq \frac{1}{c}\,,
\end{equation*}%
we obtain a contradiction, which concludes the proof.
\end{proof}

The following result is valid as $\displaystyle 
\left( \frac{n}{n+1} \right)^{n+1} < \frac{1}{e}$.
In (3.7) a non-strict inequality is also sufficient.

\vskip0.2cm\textbf{Theorem 3.4. } \textit{Assume that }$(p_{i}(n))$, $1\leq
i\leq m$, \textit{are sequences of nonnegative real numbers and} (1.1) 
\textit{holds. If}
\begin{equation}
\liminf_{n\rightarrow \infty }\sum_{i=1}^{m}\sum_{j=\tau
_{k}(n)}^{n-1}p_{i}(j) >  \frac{1}{e}\text{,} \quad, k=1, \dots, m,  \tag{3.7}
\end{equation}%
\textit{then all solutions of} (E$_{\text{R}}$) \textit{oscillate.}

\bigskip

\noindent
{\bf 3.2~~ Advanced difference equations. }
Similar oscillation theorems for the (dual) advanced difference equation (E$%
_{\text{A}}$) can be derived easily. The proofs of these theorems are
omitted, since they follow the schemes of Subsection 3.1.

\vskip0.2cm\textbf{Theorem 3.3}$^{\prime }$. \ \textit{Assume that }$%
(p_{i}(n))$, $1\leq i\leq m$, \textit{are sequences of nonnegative real
numbers, all the advances are bounded,} (1.2$^{\prime }$) \textit{holds.} 
\begin{equation}
\liminf_{n\rightarrow \infty }\sum_{i=1}^{m}\sum_{j=n+1}^{\sigma
_{k}(n)} p_{i}(j)>\left( \frac{\mu_{k} }{\mu_{k} +1}\right) ^{\mu_{k} +1}\text{,} \quad k =1, 
\dots, m,
\tag{3.8}
\end{equation}%
\textit{then all solutions of} (E$_{\text{A}}$) \textit{oscillate.}

\vskip0.2cm\textbf{Theorem 3.4}$^{\prime }$. \textbf{\ }\textit{Assume that }%
$(p_{i}(n))$, $1\leq i\leq m$, \textit{are sequences of nonnegative real
numbers, all the advances are bounded and} (1.2) \textit{holds. If}%
\begin{equation}
\liminf_{n\rightarrow \infty }\sum_{i=1}^{m}\sum_{j=n+1}^{\sigma
_{k}(n)}p_{i}(j) >  \frac{1}{e}\text{,} \quad k =1, \dots, m,  \tag{3.9}
\end{equation}%
\textit{then all solutions of} (E$_{\text{A}}$) \textit{oscillate.}


\section{EXAMPLES}

In this section, we present examples illustrating the significance of our
results. Observe that most of the relevant oscillation results cited in the
introduction cannot be applied because they assume that the delays
(advances) are constant and consequently the deviating arguments are
increasing. When possible, we compare our results to the known ones, for
variable deviations of the argument and non-monotone arguments.

\vskip0.2cm\textbf{Example 4.1. } Consider the retarded difference equation 
\begin{equation}
\Delta x(n)+px\left( \tau (n)\right) =0,\text{ \ \ \ }n\in 
\mathbb{N}
_{0}\text{,}  \tag{4.1}
\end{equation}%
where $p<1$ is a constant coefficient, and%
\begin{equation*}
\tau (n)=\left\{ 
\begin{array}{ll}
n-1=3k-1, & \text{if }n=3k \\ 
n-2=3k-1, & \text{if }n=3k+1 \\ 
n-4=3k-2, & \text{if }n=3k+2%
\end{array}%
\right. ,\text{ \ \ \ \ }k\in 
\mathbb{N}
_{0}\text{.}
\end{equation*}%
Obviously (1.1) and (2.5) are satisfied. Also, by (2.4) we have%
\begin{equation*}
\varphi (n)=\max_{0\leq s\leq n}\tau (s)=\left\{ 
\begin{array}{ll}
n-1=3k-1, & \text{if }n=3k \\ 
n-2=3k-1, & \text{if }n=3k+1 \\ 
n-3=3k-1, & \text{if }n=3k+2%
\end{array}%
\right. ,\text{ \ \ \ \ }k\in 
\mathbb{N}
_{0}\text{.}
\end{equation*}

If $n=3k$, then $\varphi (n)=n-1=$ $\tau (n)$ and, in view of (2.7), the
left-hand side in (2.10) is 
\begin{eqnarray*}
\sum_{j=\varphi (n)}^{n}p(j)a_{3}^{-1}(\varphi (n),\tau (j))
&=&\sum_{j=n-1}^{n}p(j)a_{3}^{-1}(n-1,\tau (j)) \\
&=&p\frac{1}{a_{3}(n-1,\tau (n-1))}+p\frac{1}{a_{3}(n-1,\tau (n))} \\
&=&p\frac{1}{a_{3}(n-1,n-5)}+p\frac{1}{a_{3}(n-1,n-1)} \\
&=&p\frac{1}{\mathop{\displaystyle \prod }\limits_{i=n-5}^{n-2}\left( 1-p%
\frac{1}{a_{2}(i,\tau \left( i\right) )}\right) }+p\text{,}
\end{eqnarray*}%
where 
\begin{eqnarray*}
&&\prod\limits_{i=n-5}^{n-2}\left( 1-p\frac{1}{a_{2}(i,\tau \left( i\right) )%
}\right) =\left( 1-p\frac{1}{a_{2}(n-5,\tau \left( n-5\right) )}\right)
\left( 1-p\frac{1}{a_{2}(n-4,\tau \left( n-4\right) )}\right) \times  \\
&&\times \left( 1-p\frac{1}{a_{2}(n-3,\tau \left( n-3\right) )}\right)
\left( 1-p\frac{1}{a_{2}(n-2,\tau \left( n-2\right) )}\right)  \\
&=&\left( 1-\frac{p}{a_{2}(n-5,n-7)}\right) \left( 1-\frac{p}{a_{2}(n-4,n-8)}%
\right) \left( 1-\frac{p}{a_{2}(n-3,n-4)}\right) \left( 1-\frac{p}{%
a_{2}(n-2,n-4)}\right)  \\
&=&\left( 1-\frac{p}{a_{2}(n-5,n-7)}\right) \left( 1-\frac{p}{a_{2}(n-4,n-8)}%
\right) \left( 1-\frac{p}{a_{2}(n-3,n-4)}\right) \left( 1-\frac{p}{%
a_{2}(n-2,n-4)}\right)  \\
&=&\left( 1-\frac{p}{\left( 1-\frac{p}{1-\frac{p}{\left( 1-p\right) ^{4}}}%
\right) \left( 1-\frac{p}{1-\frac{p}{1-p}}\right) }\right) ^{2}\left( 1-%
\frac{p}{1-\frac{p}{1-\frac{p}{\left( 1-p\right) ^{4}}}}\right) \times  \\
&&\times \left( 1-\frac{p}{\left( 1-\frac{p}{1-\frac{p}{\left( 1-p\right)
^{2}}}\right) ^{2}\left( 1-\frac{p}{1-\frac{p}{\left( 1-p\right) ^{4}}}%
\right) \left( 1-\frac{p}{1-\frac{p}{1-p}}\right) }\right) 
\end{eqnarray*}%
Therefore 
\begin{equation*}
\sum_{j=\varphi (n)}^{n}p(j)a_{3}^{-1}(\varphi (n),\tau (j))=p+
\end{equation*}%
\begin{equation*}
+\frac{p}{\left( 1-\frac{p}{\left( 1-\frac{p}{1-\frac{p}{\left( 1-p\right)
^{4}}}\right) \left( 1-\frac{p}{1-\frac{p}{1-p}}\right) }\right) ^{2}\left(
1-\frac{p}{\left( 1-\frac{p}{1-\frac{p}{\left( 1-p\right) ^{2}}}\right)
^{2}\left( 1-\frac{p}{1-\frac{p}{\left( 1-p\right) ^{4}}}\right) \left( 1-%
\frac{p}{1-\frac{p}{1-p}}\right) }\right) \left( 1-\frac{p}{1-\frac{p}{1-%
\frac{p}{\left( 1-p\right) ^{4}}}}\right) }.
\end{equation*}%
The computation immediately implies that, if $p\in \left(
0.174816,0.183\right) $ then%
\begin{equation*}
\limsup_{n\rightarrow \infty }\sum_{j=\varphi (n)}^{n}p(j)a_{3}^{-1}(\varphi
(n),\tau (j))>1\text{.}
\end{equation*}%
For example, for $p=0.175$ inequality (2.10) holds for $r=3$, which means
that all solutions of (4.1) oscillate. Nevertheless, for $p\in \left(
0.1741,0.183\right) $ condition (2.14) is satisfied for $r=3$, since%
\begin{equation*}
\limsup_{n\rightarrow \infty }\sum_{j=\varphi (n)}^{n}p(j)a_{3}^{-1}(\varphi
(n),\tau (j))>1-\frac{1}{2}\left( 1-p-\sqrt{1-2p-p^{2}}\right) \text{.}
\end{equation*}

Observe, however, that%
\begin{equation*}
a=\liminf_{n\rightarrow \infty }\sum_{i=\tau
(n)}^{n-1}p(i)=\liminf_{n\rightarrow \infty }\sum_{i=n-1}^{n-1}p(i)=p
\end{equation*}%
and%
\begin{equation*}
c(a)=\frac{1}{2}\left( 1-a-\sqrt{1-2a-a^{2}}\right) =\frac{1}{2}\left( 1-p-%
\sqrt{1-2p-p^{2}}\right) .
\end{equation*}%
Also,%
\begin{equation*}
\sum_{j=\varphi (n)}^{n}p(j)\prod_{i=\tau (j)}^{\varphi (n)-1}\frac{1}{1-p(i)%
}=\left\{ 
\begin{array}{ll}
p\frac{1}{\left( 1-p\right) ^{4}}+p, & \text{if }n=3k \\ 
p\frac{1}{\left( 1-p\right) ^{4}}+2p, & \text{if }n=3k+1 \\ 
p\frac{1}{\left( 1-p\right) ^{4}}+2p+p\frac{1}{1-p}, & \text{if }n=3k+2%
\end{array}%
\right. ,\text{ }k\in \mathbb{N}_{0}.
\end{equation*}%
Thus%
\begin{equation*}
\max \left\{ p\frac{1}{\left( 1-p\right) ^{4}}+p,p\frac{1}{\left( 1-p\right)
^{4}}+2p,p\frac{1}{\left( 1-p\right) ^{4}}+2p+p\frac{1}{1-p}\right\} =p\frac{%
1}{\left( 1-p\right) ^{4}}+2p+p\frac{1}{1-p}\text{.}
\end{equation*}%
If $p<0.1829$, the computation immediately implies that%
\begin{equation*}
p\frac{1}{\left( 1-p\right) ^{4}}+2p+p\frac{1}{1-p}\leq 1
\end{equation*}%
and, if $p<0.1801$, then%
\begin{equation*}
p\frac{1}{\left( 1-p\right) ^{4}}+2p+p\frac{1}{1-p}\leq 1-\frac{1}{2}\left(
1-p-\sqrt{1-2p-p^{2}}\right) \text{.}
\end{equation*}%
Therefore, if $p\in \left( 0.174816,0.1801\right) $, the conditions (2.1)
and (2.2) are not satisfied.

Moreover, observe that if $p\in \left( 0.1741,0.174816\right) $, condition
(2.14) holds for $r=3$ while (2.10 ) and also (2.1) and (2.2) are not
satisfied. For example, for $p=0.1742$ only (2.14) holds for $r=3$.

\vskip0.2cmAt this point we should remark that in condition (2.14) an extra
requirement is needed: $0<\alpha \leq 1/e.$ Thus, when $\alpha \rightarrow 0$
(2.14) coincides with (2.10).

\vskip0.2cm\textbf{Example 4.2. } Consider the delay difference equation 
\begin{equation}
\Delta x(n)+\frac{1}{8}x\left( \tau _{1}(n)\right) +\frac{1}{12}x\left( \tau_{2}(n)\right) =0,
\text{ \ \ \ }n\geq 0\text{,}  \tag{4.2}
\end{equation}
with
\begin{equation*}
\tau _{1}(n)=\left\{ 
\begin{array}{ll}
n-3 & \text{if }n\text{ is even} \\ 
n-1 & \text{if }n\text{ is odd}
\end{array}
\right. \text{ \ \ \ \ and \ \ \ \ }\tau _{2}(n)=\left\{ 
\begin{array}{ll}
n-4 & \text{if }n\text{ is even} \\ 
n-1 & \text{if }n\text{ is odd}%
\end{array}%
\right. \text{.}
\end{equation*}
Here, it is clear that (1.1) and (2.5) are satisfied. Also, by (2.3) and
(2.4) we have
\begin{equation*}
\varphi _{1}(n)=\varphi _{2}(n)=\left\{ 
\begin{array}{ll}
n-2 & \text{if }n\text{ is even} \\ 
n-1 & \text{if }n\text{ is odd}
\end{array}
\right.
\end{equation*}
and
\begin{equation*}
\varphi (n)=\max_{1\leq i\leq 2}\varphi _{i}(n)=\left\{ 
\begin{array}{ll}
n-2 & \text{if }n\text{ is even} \\ 
n-1 & \text{if }n\text{ is odd}
\end{array}
\right. \text{.}
\end{equation*}
If $n$ is even, then $\varphi (n)=n-2$, $\tau _{1}(n)=n-3$, $\tau
_{2}(n)=n-4 $ and, in view of (2.7), (2.10) gives%
\begin{equation*}
\sum_{j=\varphi (n)}^{n}\sum_{i=1}^{m}p_{i}(j)a_{1}^{-1}(\varphi (n),\tau
_{i}(j))=\sum_{j=n-2}^{n}\sum_{i=1}^{2}p_{i}(j)a_{1}^{-1}(n-2,\tau _{i}(j))
\end{equation*}%
\begin{equation*}
=\frac{1}{8}\cdot \left[ \frac{1}{\left( 1-\frac{5}{24}\right) ^{3}}+1+\frac{%
1}{1-\frac{5}{24}}\right] +\frac{1}{12}\cdot \left[ \frac{1}{\left( 1-\frac{5%
}{24}\right) ^{4}}+1+\frac{1}{\left( 1-\frac{5}{24}\right) ^{2}}\right]
\simeq 0.963276895\text{.}
\end{equation*}%
Also%
\begin{equation*}
\sum_{j=\varphi (n)}^{n}\sum_{i=1}^{m}p_{i}(j)a_{2}^{-1}(\varphi (n),\tau
_{i}(j))=\sum_{j=n-2}^{n}\sum_{i=1}^{2}p_{i}(j)a_{2}^{-1}(n-2,\tau _{i}(j))
\end{equation*}%
\begin{eqnarray*}
&=&\frac{1/8}{\left[ 1-\frac{1}{8}\cdot \frac{1}{1-\frac{5}{24}}-\frac{1}{12}%
\cdot \frac{1}{1-\frac{5}{24}}\right] ^{2}\left[ 1-\frac{1}{8}\cdot \frac{1}{%
\left( 1-\frac{5}{24}\right) ^{3}}-\frac{1}{12}\cdot \frac{1}{\left( 1-\frac{%
5}{24}\right) ^{4}}\right] } \\
&&+\frac{1/12}{\left[ 1-\frac{1}{8}\cdot \frac{1}{\left( 1-\frac{5}{24}%
\right) ^{3}}-\frac{1}{12}\cdot \frac{1}{\left( 1-\frac{5}{24}\right) ^{4}}%
\right] ^{2}\left[ 1-\frac{1}{8}\cdot \frac{1}{1-\frac{5}{24}}-\frac{1}{12}%
\cdot \frac{1}{1-\frac{5}{24}}\right] ^{2}} \\
&&+\frac{1}{8}+\frac{1}{12}+\frac{1/8}{\left[ 1-\frac{1}{8}\cdot \frac{1}{1-%
\frac{5}{24}}-\frac{1}{12}\cdot \frac{1}{1-\frac{5}{24}}\right] } \\
&&+\frac{1/12}{\left[ 1-\frac{1}{8}\cdot \frac{1}{\left( 1-\frac{5}{24}%
\right) ^{3}}-\frac{1}{12}\cdot \frac{1}{\left( 1-\frac{5}{24}\right) ^{4}}%
\right] \left[ 1-\frac{1}{8}\cdot \frac{1}{1-\frac{5}{24}}-\frac{1}{12}\cdot 
\frac{1}{1-\frac{5}{24}}\right] }
\end{eqnarray*}
i.e.,
\begin{equation*}
\sum_{j=\varphi (n)}^{n}\sum_{i=1}^{m}p_{i}(j)a_{2}^{-1}(\varphi (n),\tau
_{i}(j))\simeq 1.553022949
\end{equation*}%
Therefore 
\begin{equation*}
\limsup_{n\rightarrow \infty }\sum_{j=\varphi
(n)}^{n}\sum_{i=1}^{m}p_{i}(j)a_{2}^{-1}(\varphi (n),\tau _{i}(j))\geq
1.553022949>1\text{,}
\end{equation*}
that is, condition (2.10) of Theorem 2.4 is satisfied and therefore all
solutions of (4.2) oscillate. Also,%
\begin{equation*}
a_{1}=\liminf_{n\rightarrow \infty }\sum_{i=\varphi _{1}(n)}^{n-1}p_{1}(i)=%
\frac{1}{8}\text{ \ \ \ \ and \ \ \ \ }a_{2}=\liminf_{n\rightarrow \infty
}\sum_{i=\varphi _{2}(n)}^{n-1}p_{2}(i)=\frac{1}{12}
\end{equation*}
and therefore
\begin{equation*}
\alpha =\min \left\{ \alpha _{i}:1\leq i\leq 2\right\} =\frac{1}{12}.
\end{equation*}%
Hence%
\begin{equation*}
\limsup_{n\rightarrow \infty }\sum_{j=\varphi
(n)}^{n}\sum_{i=1}^{m}p_{i}(j)a_{2}^{-1}(\varphi (n),\tau _{i}(j))>1-\frac{1%
}{2}\left( 1-a-\sqrt{1-2a-a^{2}}\right) \simeq 0.996196338\text{,}
\end{equation*}%
that is, condition (2.14) of Theorem 2.5 is also satisfied.

Here, it is obvious that%
\begin{equation*}
n-\tau _{1}(n)\leq 3=M_{1}\text{ \ \ \ \ and \ \ \ \ }n-\tau _{2}(n)\leq
4=M_{2}
\end{equation*}%
and therefore%
\begin{equation*}
M=\max_{1\leq i\leq 2}M_{i}=4\text{.}
\end{equation*}%
Observe that%
\begin{equation*}
\liminf_{n\rightarrow \infty }\sum_{i=1}^{2}\sum_{j=\tau
_{i}(n)}^{n-1}p_{i}(j)=\frac{1}{8}+\frac{1}{12}=\frac{5}{24}<\left( \frac{M}{%
M+1}\right) ^{M+1}=0.32768\text{.}
\end{equation*}%
That is, condition (3.3) is not satisfied. Also the second condition in
(3.1) is not satisfied.

Fig.~\ref{figure1} below illustrates a couple of solutions of (4.2). As can be
observed, both solutions are oscillatory.


\begin{figure}[ht]
\centering
\includegraphics[width=0.3 \linewidth]{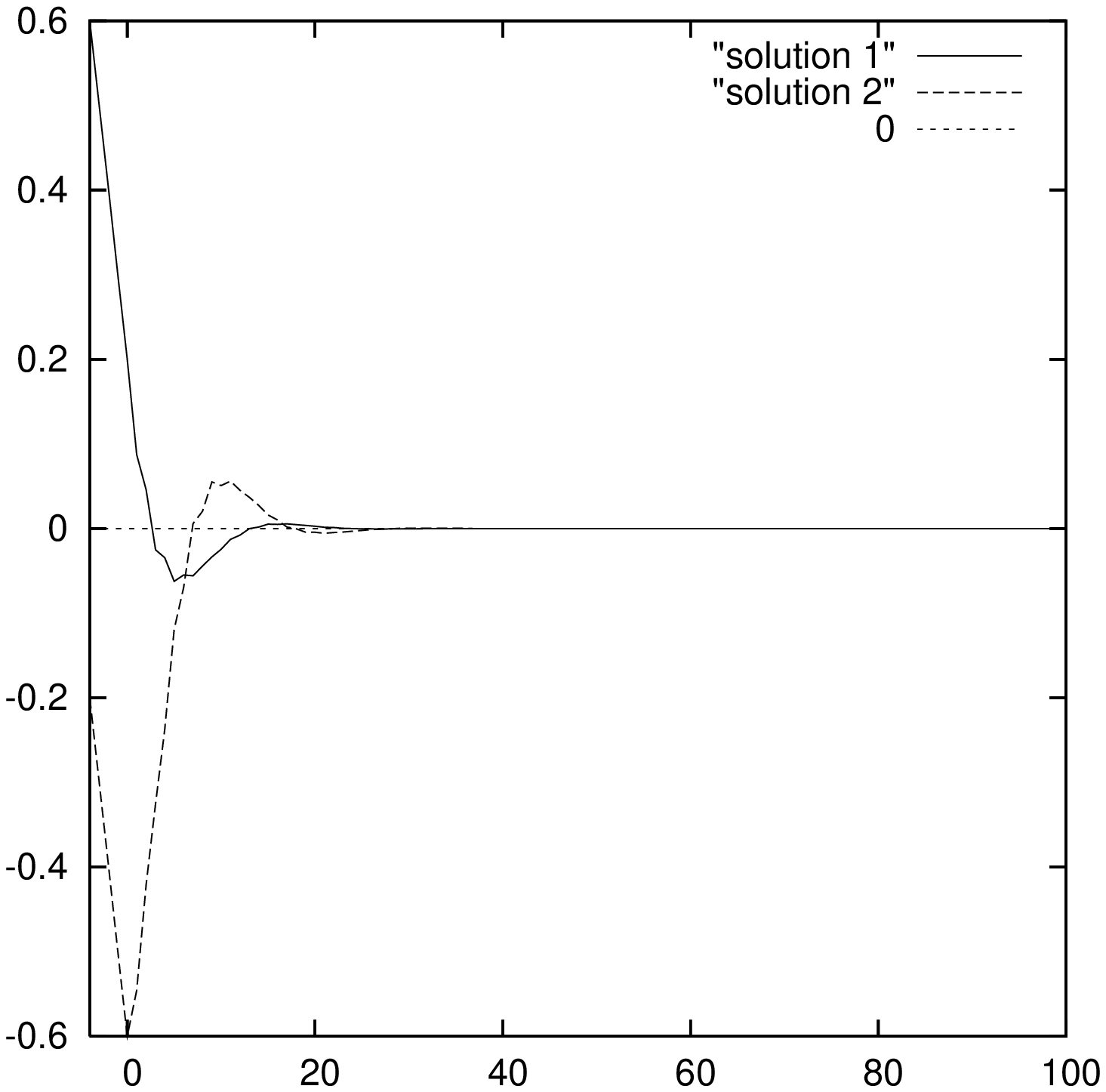}~~~~~
\includegraphics[width=0.3 \linewidth]{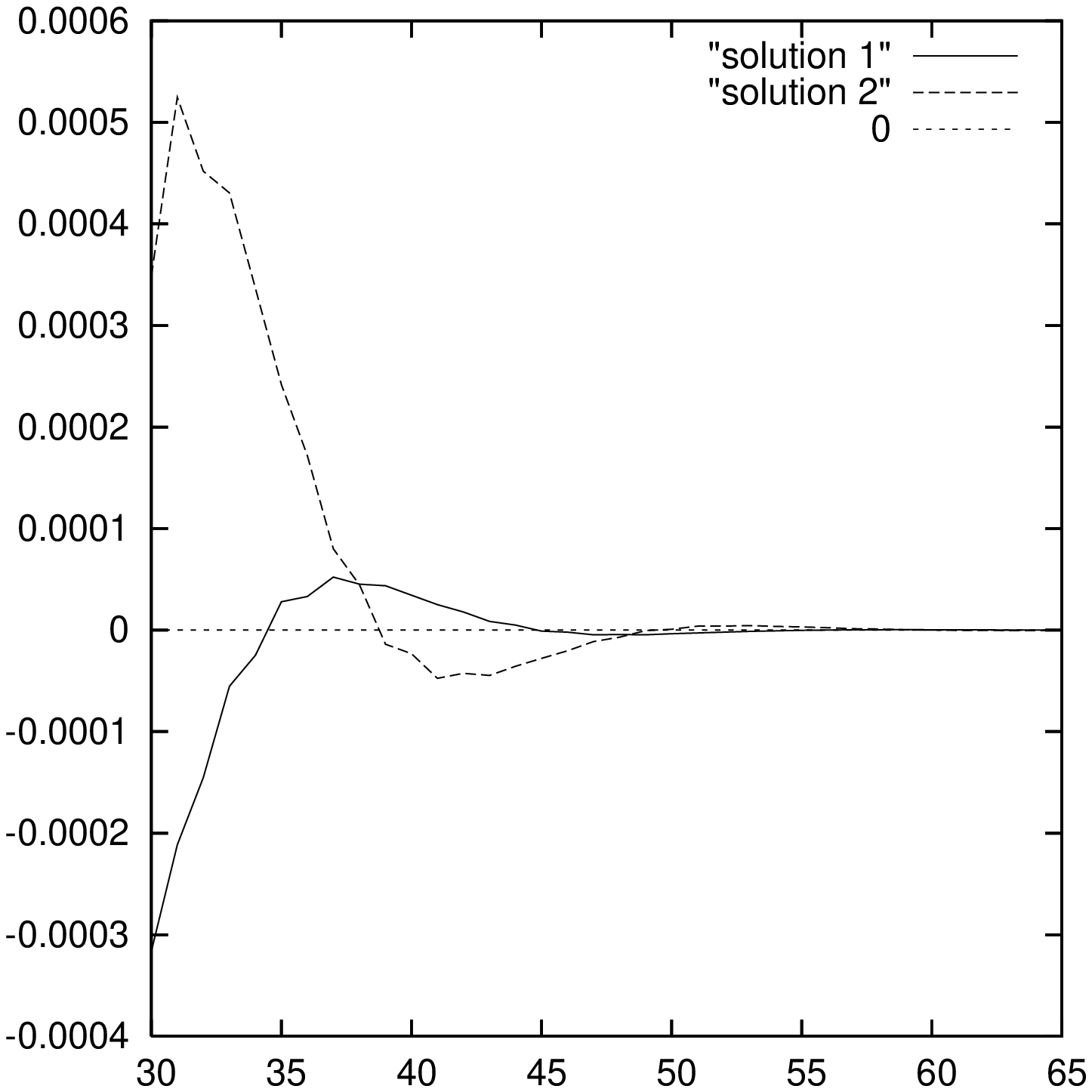}~~~~~
\includegraphics[width=0.3 \linewidth]{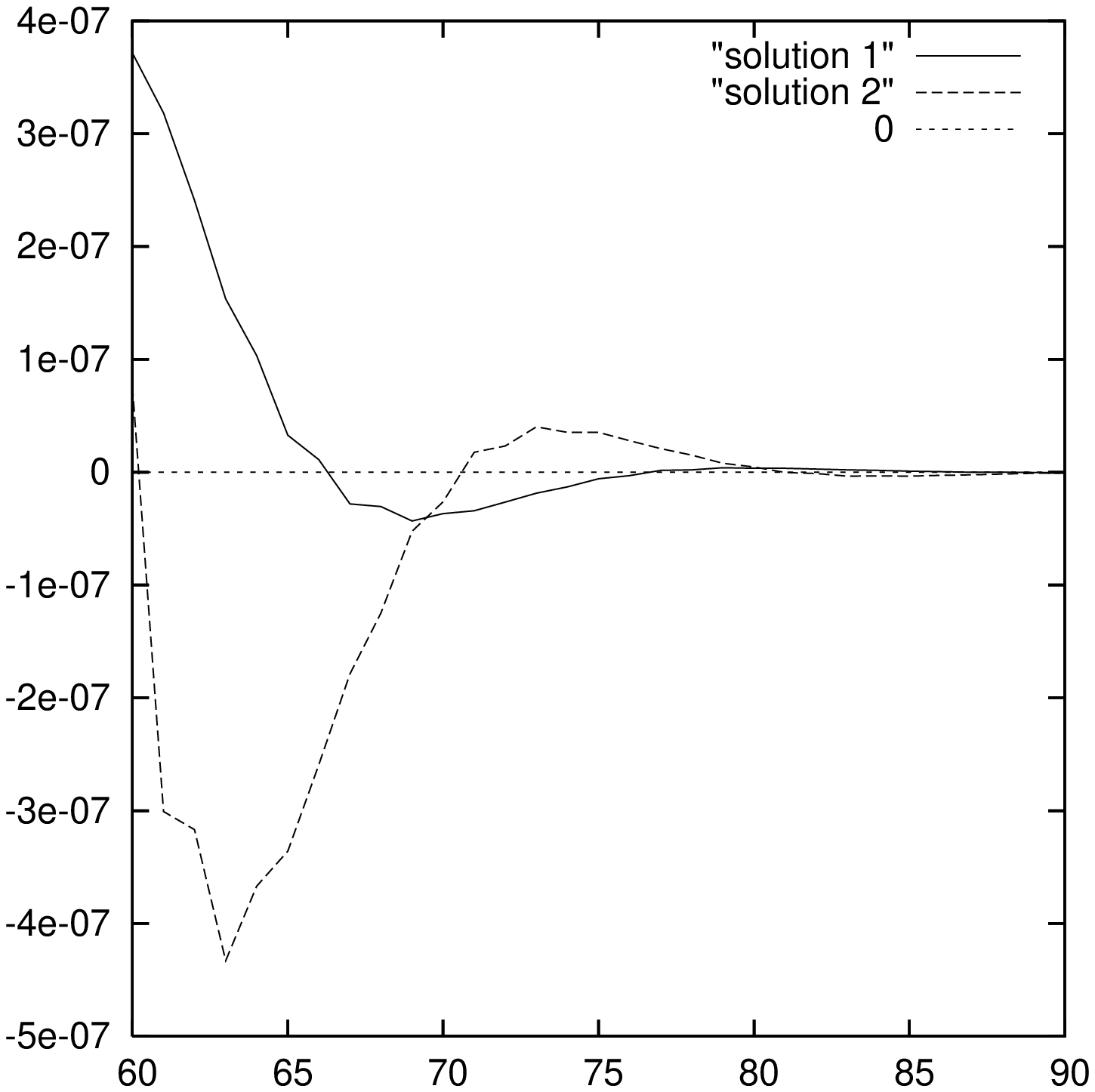}
\caption{Two solutions of (4.2) for (left) $-4 \leq n \leq 100$, (middle)
$30 \leq n \leq 65$, (right) $60  \leq n \leq 90$. Both solutions are oscillatory.} 
\label{figure1}
\end{figure}


\vskip0.2cm\textbf{Example 4.3. } 
Consider the delay difference equation 
\begin{equation}
\Delta x(n)+ a(n) x\left( \tau _{1}(n)\right) +\frac{1}{125}x\left(\tau _{2}(n)\right) =0,\text{ \ \ \ }n\geq 0\text{,}  \tag{4.3}
\end{equation}%
with
\begin{equation*}
\tau _{1}(n)=\left\{ 
\begin{array}{ll}
n-1 & \text{if }n\text{ is even,} \\ 
n-2 & \text{if }n\text{ is odd,}
\end{array}%
\right. \ \  \tau _{2}(n)=\left\{ 
\begin{array}{ll}
n-2 & \text{if }n\text{ is even,} \\ 
n-3 & \text{if }n\text{ is odd,}
\end{array}
\right. 
\end{equation*}
\begin{equation*}
a(n)=\left\{
\begin{array}{ll}   
\frac{3}{125} & \text{if }n\text{ is even,} \\ \\
\frac{37}{125} & \text{if }n\text{ is odd.}   
\end{array}%
\right.
\end{equation*}
Evidently
\begin{equation*}
n-\tau _{1}(n)\leq 2=M_{1}\text{ \ \ \ \ and \ \ \ \ }n-\tau _{2}(n)\leq 3=M_{2}
\end{equation*}%
and therefore for $k=1$
$$
\liminf_{n\rightarrow \infty }\sum_{i=1}^{2}\sum_{j=\tau_{1}(n)}^{n-1}p_{i}(j)
=\frac{37}{125}+ \frac{1}{125}=\frac{38}{125}=0.304 >\left( \frac{M_1}{M_1+1}\right) ^{M_1+1} \approx 0.2962963,
$$
while for $k=2$
$$
\liminf_{n\rightarrow \infty }\sum_{i=1}^{2}\sum_{j=\tau_{2}(n)}^{n-1}p_{i}(j)
= \frac{37}{125}+ \frac{3}{125} + \frac{2}{125} = 0.336
>\left( \frac{M_2}{M_2+1}\right)^{M_2+1} \approx 0.31640625\text{,}
$$
that is, condition (3.3) is satisfied and, by Theorem 3.3, all solutions of
Eq.(4.3) oscillate.


\medskip

\centerline{\bf Acknowledgments}

\medskip

The first author was partially supported by NSERC, grants RGPIN/261351-2010 and RGPIN-2015-05976.

\end{document}